# Non-limit integration of differential equations. General solution for van der Pol equation

Sergey V. Zuev[1]

## Introduction

Using the physical background based on the earlier author's researches the new method of the ODE and PDE integration is proposed. As illustrations of the method the general solutions for free oscillator and for heat equation are derived. As a result of the article the general solution for van der Pol equation is constructed.

## Inertial frames and differential calculus

Let us consider two inertial reference frames (IRF) with coordinates (include time) of the same test event as $(x_1, y_1, z_1, t_1)$, $(x_2, y_2, z_2, t_2)$. We assume that general relativity theory (GRT) is valid and all the intervals between the three events: the reference point in the IRF-1, the reference point in the IRF-2 and the test event does not depend on the choice of IRF. It is easy to show that in this case the initial event in the IRF-2 with respect to the IRF-1 corresponds to the coordinates $(x^1 = x = x_1 - x_2, x^2 = y = y_1 - y_2, x^3 = z = z_1 - z_2, x^4/c = t = t_1 - t_2)$ with respect of any kind of space-time metric.

From physical considerations it is clear that the bodies of reference of both IRF have a nonzero rest mass. That is why they move relative to each other with speed strictly less than $c$. In case of the absence of the second IRF, the motion of a test particle in the first IRF would be described in terms of the SRT, which corresponds to the flat Minkowski space-time metric. Hence, the required space of possible IRF-2 is to be asymptotically Minkowski space-time. The group of isometries for the Minkowski space is Poincare group. While noting the position of the origin of coordinates and time for the IRF-1, we have the possibility of translation in the target space exhausted. That is why the isometry group is reduced to Lorentz group. Because of the principle of causality it is awaited that space-time is oriented: causality yields space-time to be oriented in time and the interval invariance leads to the whole space-time to be oriented. All that require the isometry group to be $SO(3,1)$.

For oriented pseudo-Riemannian spaces with isometry group $SO(3,1)$ the next proposition is valid.

**Theorem 1** *[1,2]. Let U be a simply connected oriented area in some pseudo-Riemannian 4-manifold with isometry group $SO(3,1)$. Let $(x^1, x^2, x^3, x^4)$ be some coordinates on U such that, the metric $g$ of the manifold tends to the Minkowski metric when $r \equiv {x^4}^2 - {x^3}^2 - {x^2}^2 - {x^1}^2 \to \infty$. Then*
$$g_{ij} = u'' x^k x^l \eta_{ik} \eta_{jl} + u' \eta_{ij}$$
*where*
$$u' = \frac{(r^4 - a^4)^{1/4}}{r}, u'' = \frac{du'}{dr}, \eta_{ij} = diag(-1,-1,-1,1), a = const \neq 0$$
*Inverse, if the metric tensor of the pseudo-Riemannian 4-manifold with signature $(---+)$ in every chart U could take an origin of $g_{ij}$, than the manifold has isometry group $SO(3,1)$.*

[1] Belgorod State Technological University, Belgorod, Russia

Therefore, when one takes time-like area of space-time, there exists some (positive) parameter $a$, which limits minimum of $r$. In case $r = a$ the space-time metrics and its invariants become infinite, i.e. there is singularity. For values $r < a$ the interval between reference events in IRF-1 and IRF-2 becomes space-like and that is why there is no any correlations between IRF-1 and IRF-2.

The circumstance that $r$ is bounded from below does not prohibit for IRF-1 and IRF-2 to have common reference body. In that case $r = x^{4^2} = c^2(t_1 - t_2)^2$, i.e. for such pair of IRFs there exists a the minimum time interim that can separate any of two events with respect of those IRFs. The appointed interim is equal to $\frac{\sqrt{a}}{c}$ and significantly non-zero.

For any reference body (and for any object with non-zero rest mass) the all range of parameters independent from coordinates and self-time is exhausted by rest mass and gauge invariants only. Concerning gauge invariants one can say that they are not directly connected with space-time. Therefore in case of IRFs connected to two different reference bodies one has $a = F(m_1, m_2)$, where $F$ is some function. When two IRFs connected to the only reference body one has $a = f(m)$.

It is important that for any object in any IRF with reference body rest mass $m$ there exists the minimal interim which separates any event with the object and in the same time "registrable" with respect of the IRF. The interim depends on $m$ only and it is strong positive.

While we remember that in classical differential calculus in physics one has no any minimum limit for time: any derivative is calculated in limit $t \to 0$ (the quantum mechanics is looking like an exemption but we dealing with classical case now). That is why the follow proposition proven:

*The differential calculus in physics (using time- and space- derivatives) is not exact description for physical processes in IRFs. The differential calculus gives approximation accuracy only.*

Hence it is possible to find some effects in processes which could not be described by differential equations with ordinary derivatives.

**The method of non-limit integration for differential equations**

Instead of classical definition of derivative $\dot{x}(t)$ of function $x(t)$ let us consider the following:
$$\dot{x}(t) = \frac{x(t + \tau) - x(t)}{\tau}$$
where $\tau$ is some non-zero positive minimal possible value of $t$. Unlike to ordinary derivative, let us call the function $\dot{x}(t)$ *non-limit derivative* of $x(t)$. One should especially appoint the basic properties of non-limit derivative because of some of them are quite different from the traditional ones.
- Non-limit derivative of constant is zero, non-limit derivative of identity is equal to one:
$$\dot{C} = 0 \qquad \dot{t} = 1$$
- Non-limit derivative of linear combination of functions is linear combination of non-limit derivative of these functions:
$$\dot{(kf(t) + l\, g(t))} = k\dot{f}(t) + l\dot{g}(t)$$
- Non-limit derivative of the functions product:
$$\dot{g(t)f(t)} = \dot{g}(t)f(t) + g(t)\dot{f}(t) + \tau\dot{g}(t)\dot{f}(t)$$
- Non-limit derivative of the functions quotient:

$$\left(\frac{f(t)}{g(t)}\right)^{\cdot} = \frac{g(t)\dot{f}(t) - f(t)\dot{g}(t)}{g(t)(g(t) + \tau \dot{g}(t))}$$

- Some examples of non-limit derivative of elementary functions:
$$\dot{p^x} = \frac{p^x(p^{x\tau} - 1)}{\tau}, \quad \forall p$$
$$\dot{x^n} = x(\dot{x^{n-1}}) + \dot{x}(\tau\dot{x} + x)^{n-1}$$

(just recurrently)
$$(\log_k x)^{\cdot} = \frac{1}{\tau}\log_k\left(1 + \tau\frac{\dot{x}}{x}\right)$$

*Non-limit differential* is finite and determined by
$$\Delta x = \tau \dot{x}$$

Let us call any equation containing non-limit derivatives instead ordinary derivatives as *non-limit differential equation* (NLDE). When deal with non-limit derivatives with respect of one variable the NLDE will be NLODE and otherwise become NLPDE. Both types are exampled below. The simplest method for NLDE integration is change of variable. The method is illustrated by two samples below. The main goal of the method is reduce the equation to an identity using change(s) of dependent or independent variable(s). Another one method is illustrated by van der Pol equation integration. Together with changes of variable one used analysis of sequence which is approximated by unknown function.

The transition from NLDE solution(s) to the ordinary ODE or PDE solution(s) one may performs by the limit $\tau \to 0$ or by $\tau$ exemption using suitable change of variables.

### Differential equations for physical systems

Let us consider differential equations which describe the following 1-dimensional (for simplicity) physical systems: free oscillations and heat transfer.

**Free oscillations.** The system is usually described by the next ordinary differential equation of second order:
$$\ddot{x} + \omega^2 x = 0$$
where $\omega$ is some positive constant with physical sense as frequency. It is easy to find the general solution for this equation:
$$x = Ae^{i\omega t} + Be^{-i\omega t}$$
or
$$x = C_1 \sin \omega t + C_2 \cos \omega t$$
for real values.

In case of discrete time the system might be described using the equation of the same view but with the non-limit derivatives instead. Let us derive the general solution of the equation for sample. While the only elementary function which has a proportional function as a non-limit derivative is exponential function $p^x$, it is natural to try to find the solution looks like $x = Ap^{Bt}$. Using the variable $x$ change:
$$x = Ap^{Bt}$$
(A and B are some constants) we find
$$\dot{x} = Ap^{Bt}\frac{1}{\tau}(p^{B\tau} - 1)$$
and
$$\ddot{x} = Ap^{Bt}\frac{1}{\tau^2}(p^{B\tau} - 1)^2$$

Substituting into the original equation we get (if $A \neq 0$)
$$(p^{B\tau} - 1)^2 = -\omega^2 \tau^2$$
hence
$$p^{B\tau} = 1 \pm i\tau\omega$$
$$B = \frac{1}{\tau}\log_p(1 \pm i\tau\omega)$$
and because of the linearity of the equation the general solution is a linear combination of exponential functions of the following form:
$$x = C_1(1 + i\tau\omega)^{\frac{t}{\tau}} + C_2(1 - i\tau\omega)^{\frac{t}{\tau}}$$
It is evident that when $\tau \to 0$ the appointed $x$ goes to the well-known general solution of exponential form.

**The heat equation.** Now let us consider heat equation for homogeneous rod without any heat sources:
$$\frac{\partial u}{\partial t} - \alpha \frac{\partial^2 u}{\partial x^2} = 0$$

It's general solution was found in [3] in form of proved convergent sum:
$$u(x,t) = \sum_{n=1}^{\infty} \left( \frac{2}{l} \int_0^l \varphi(\xi) \sin\left(\frac{\pi n}{l}\xi\right) d\xi \right) \sin\left(\frac{\pi n}{l}x\right) \exp\left(-a^2 \left(\frac{\pi n}{l}\right)^2 t\right)$$

As it showed below, this solution is a partial case of another one. But it is not shown that the new solution is not equivalent to that.

In general the steps of non-limit differential are different for $x$ and $t$. Let us denote as $\xi$ a step on $x$ and as $\tau$ a step on $t$. Using variable change $y = \frac{\tau}{\xi}x$, we find the equation
$$\frac{\partial u}{\partial t} - \beta \frac{\partial^2 u}{\partial y^2} = 0$$
where $\beta = \alpha \frac{\tau^2}{\xi^2}$ and the non-limit differential steps for both independent variables becomes equal to $\tau$. This equation is linear and the linear combination of its solutions is a solution. Now let us note that ranges of values for variables with the same differential step coincide with each other and that is why the range of target function is completely determined by that range. Appointed circumstance allows us to try to find the solution in form of:
$$u(t,y) = C_1(\tau\lambda_1 + 1)^{t/\tau}(\tau\lambda_1 + 1)^{y/\tau} + C_2(\tau\lambda_1 + 1)^{t/\tau}(\tau\lambda_2 + 1)^{y/\tau}$$
$$+ C_3(\tau\lambda_2 + 1)^{t/\tau}(\tau\lambda_1 + 1)^{y/\tau} + C_4(\tau\lambda_2 + 1)^{t/\tau}(\tau\lambda_2 + 1)^{y/\tau}$$
where $\lambda_1$ and $\lambda_2$ are solutions of the characteristic equation
$$\lambda - \beta\lambda^2 = 0$$
so that
$$\lambda_1 = 0, \quad \lambda_2 = \frac{1}{\beta}$$
and
$$u(t,y) = C_1 + C_2\left(\frac{\tau}{\beta} + 1\right)^{y/\tau} + C_3\left(\frac{\tau}{\beta} + 1\right)^{t/\tau} + C_4\left(\frac{\tau}{\beta} + 1\right)^{t/\tau}\left(\frac{\tau}{\beta} + 1\right)^{y/\tau}$$
Substitution into the equation gives
$$C_2 = C_3 = 0$$
and $C_1$, $C_4$ can be an arbitrary constants. Finally:
$$u(t,y) = C_1 + C_4\left(\frac{\tau}{\beta} + 1\right)^{t/\tau}\left(\frac{\tau}{\beta} + 1\right)^{y/\tau}$$
Recall that $\beta = \alpha \frac{\tau^2}{\xi^2}$ and $y = \frac{\tau}{\xi}x$ that leads to

$$u(t,x) = C_1 + C_4\left(\frac{\xi^2}{\alpha\tau} + 1\right)^{\frac{t}{\tau}+\frac{x}{\xi}}$$

This is exact solution of heat equation for homogeneous rod without heat sources. However $\tau$ and $\xi$ are still undetermined. To get the limited (traditional) solution, it is necessary to put $\tau, \xi \to 0$. Let us denote $\lim_{\tau,\xi\to 0} \frac{\xi}{\tau} = k$ in the same time. Then one has:

$$u(t,x) = C_1 + C_4 e^{\frac{k^2 t}{\alpha}} e^{\frac{kx}{\alpha}} \qquad (1)$$

The known general solution of the heat equation is determined for $x \in (0; l), \tilde{u}(0,x) = \varphi(x)$ and $\tilde{u}(t,0) = \tilde{u}(t,l) = 0$ (see [3]):

$$\tilde{u}(t,x) = \sum_{n=1}^{\infty} \left(\frac{2}{l}\int_0^l \varphi(\xi)\sin\left(\frac{\pi n}{l}\xi\right)d\xi\right)\sin\left(\frac{\pi n}{l}\xi\right)\sin\left(\frac{\pi n}{l}x\right)\exp\left(-\alpha\left(\frac{\pi n}{l}\right)^2 t\right).$$

If we put in (1) $k_n = i\alpha\frac{\pi n}{l}$ then this solution will be the imaginary part of the sum of the solutions of (1) type:

$$\tilde{u}(t,x) = \sum_{n=1}^{\infty} A_n e^{\frac{k_n^2 t}{\alpha}} e^{\frac{k_n x}{\alpha}}$$

where $A_n = \frac{2}{l}\int_0^l \varphi(\xi)\sin\left(\frac{\pi n}{l}\xi\right)d\xi$. But one could construct also the real part of the sum and the result will also be a solution. But it is not evident that the appointed solutions are not equivalent.

**Integration of van der Pol equation for arbitrary values of parameters**

The classical van der Pol equation [7] describes stable relaxation oscillations (otherwise known as *limit cycles*) and has the following form

$$\ddot{x} - \lambda(1-x^2)\dot{x} + \omega^2 x = 0$$

This equation has physical applications in such fields as neuron networks, determined chaos theory and some other aspects of nonlinear dynamics. The van der Pol equation is a very good example to show the nature of non-linearity, that is why much attention for the equation paid in monographs and reviews on nonlinear dynamics (see, for instance [4],[5], and [6] as well as references therein). The general solution of the equation for arbitrary $\lambda$ and $\omega$ values not found yet. Let us find it.

First of all, let us consider the van der Pol equation as NLDE. Let us change the depended variable:

$$y = \log_k x \quad \text{or} \quad x = k^y$$

where $k$ is still arbitrary constant. Dividing by $x$, one has

$$\frac{\ddot{x}}{x} - \lambda(1-x^2)\frac{\dot{x}}{x} + \omega^2 = 0$$

note that

$$\frac{\dot{x}}{x} = \frac{k^{\dot{y}\tau} - 1}{\tau} \quad \text{and}$$

$$\frac{\ddot{x}}{x} = \frac{k^{\dot{y}\tau}\left(k^{\ddot{y}\tau^2}k^{\dot{y}\tau} - 2\right) + 1}{\tau^2}$$

Substituting into the equation, multiplying by $\tau^2$ and giving similar, we find

$$k^{\ddot{y}\tau^2 + 2\dot{y}\tau} - (2 + \lambda\tau)k^{\dot{y}\tau} + \lambda\tau k^{\dot{y}\tau + 2y} - \lambda\tau k^{2y} + 1 + \lambda\tau + \tau^2\omega^2 = 0. \qquad (*)$$

Note that it is linear combination of exponential functions which are known to be linearly independent functions. Hence some factors must be proportional with constant coefficients and at least one of them must be independent from $t$. It is clear that function $y$ not permitted to be constant. Therefore while

$$\Lambda \equiv 1 + \lambda\tau + \tau^2\omega^2 \neq 0,$$

only three following expressions might be equal to constants

$$\dot{y}\tau$$
$$\ddot{y}\tau^2 + 2\dot{y}\tau$$

or

$$2y + \dot{y}\tau.$$

Let us consider all cases.

a. Let $\dot{y}\tau = const = a_2$. Then $y = \frac{a_2}{\tau}t + b_2$ and (*) becomes

$$k^{2a_2} - (2+\lambda\tau)k^{a_2} + \lambda\tau k^{a_2+2\frac{a_2t}{\tau}+2b_2} - \lambda\tau k^{2\frac{a_2t}{\tau}+2b_2} + 1 + \lambda\tau + \tau^2\omega^2 = 0.$$

It is evident that in that case $a_2 = 0$ and $y = const$. So this case does not lead to the solution.

b. Let us consider case $\ddot{y}\tau^2 + 2\dot{y}\tau = a_1 = const = a_1$. Then it is easy to verify that

$$y = c_1(-1)^{\frac{t}{\tau}} + \frac{a_1}{2\tau}t + \frac{b_1 - \frac{a_1}{2}}{2}$$ and equation (*) takes the following form:

$$k^{a_1} - (2+\lambda\tau)k^{\frac{a_1}{2}-2c_1(-1)^{\frac{t}{\tau}}} + \lambda\tau k^{\frac{a_1}{\tau}t+b_1} - \lambda\tau k^{\frac{a_1t}{\tau}+2c_1(-1)^{\frac{t}{\tau}}+b_1-\frac{a_1}{2}} + 1 + \lambda\tau + \tau^2\omega^2 = 0$$

While $t$ multiples $\tau$, just two cases permitted for real values: $(-1)^{\frac{t}{\tau}} = +1$ and $(-1)^{\frac{t}{\tau}} = -1$. For $(-1)^{\frac{t}{\tau}} = +1$ we have

$$k^{a_1} - (2+\lambda\tau)k^{\frac{a_1}{2}-2c_1} + \lambda\tau k^{\frac{a_1}{\tau}t}\left(k^{b_1} - k^{2c_1+b_1-\frac{a_1}{2}}\right) + 1 + \lambda\tau + \tau^2\omega^2 = 0$$

and hence

$$2c_1 - \frac{a_1}{2} = 0.$$

But in case of $(-1)^{\frac{t}{\tau}} = -1$ we derive from (*)

$$k^{a_1} - (2+\lambda\tau)k^{\frac{a_1}{2}+2c_1} + \lambda\tau k^{\frac{a_1}{\tau}t}\left(k^{b_1} - k^{-2c_1+b_1-\frac{a_1}{2}}\right) + 1 + \lambda\tau + \tau^2\omega^2 = 0$$

and $2c_1 + \frac{a_1}{2} = 0$. Of course, there is no any possibility to be $2c_1 - \frac{a_1}{2} = 0$ and $2c_1 + \frac{a_1}{2} = 0$ simultaneously if $a_1$ and $c_1$ are both non-zero. So this case is also trivial.

c. Let us consider finally $\dot{y}\tau + 2y = a_3$, that leads to

$$y = \frac{a_3}{2} + b_3(-1)^{\frac{t}{\tau}}$$
$$\dot{y} = -\frac{2b_3}{\tau}(-1)^{\frac{t}{\tau}}$$

The (*) equation becomes

$$-(2+\lambda\tau)k^{-2b_3(-1)^{\frac{t}{\tau}}} + \lambda\tau k^{a_3} - \lambda\tau k^{a_3+2b_3(-1)^{\frac{t}{\tau}}} + 2 + \lambda\tau + \tau^2\omega^2 = 0$$

and could be divided onto two equations:

$$-(2+\lambda\tau)k^{2b_3} + \lambda\tau k^{a_3} - \lambda\tau k^{a_3-2b_3} + 2 + \lambda\tau + \tau^2\omega^2 = 0 \quad \text{when} \quad (-1)^{\frac{t}{\tau}} = -1$$

and

$$-(2+\lambda\tau)k^{-2b_3} + \lambda\tau k^{a_3} - \lambda\tau k^{a_3+2b_3} + 2 + \lambda\tau + \tau^2\omega^2 = 0 \quad \text{when} \quad (-1)^{\frac{t}{\tau}} = +1$$

The solution will be generated by such pair of $a_3, b_3$ which pay in the identity both of this two equations (i.e. it needs to solve the system of two algebraic equations).
Solving the system, one has:

$$x = k^y = \sqrt{k^{a_3}}\left(\sqrt{k^{2b_3}}\right)^{(-1)^{\frac{t}{\tau}}}$$

where

$$k^{a_3} = \frac{2+\lambda\tau}{\lambda\tau}, \quad k^{2b_3} = 1 + \frac{\omega^2\tau^2}{2(2+\lambda\tau)} \pm \sqrt{\left(1+\frac{\omega^2\tau^2}{2(2+\lambda\tau)}\right)^2 - 1}$$

and we can derive the solution of the van der Pol equation in the following brief record

$$x = \sqrt{\frac{2+\lambda\tau}{\lambda\tau}} \left[\Omega \pm \sqrt{\Omega^2 - 1}\right]^{\frac{(-1)^{t/\tau}}{2}} \quad (**)$$

where $\Omega = 1 + \frac{\omega^2 \tau^2}{2(2+\lambda\tau)}$ and $\tau \neq \frac{-\lambda \pm \sqrt{\lambda^2 - 4\omega^2}}{2\omega^2}$ because of $\Lambda \neq 0$.

In case of is $\Lambda \neq 0$ the (*) equation becomes to be
$$k^{\ddot{y}\tau^2 + 2\dot{y}\tau} - (2+\lambda\tau)k^{\dot{y}\tau} + \lambda\tau k^{\dot{y}\tau + 2y} - \lambda\tau k^{2y} = 0$$
and there are 6 conditions for $k$ factors:

a. $\ddot{y}\tau^2 + 2\dot{y}\tau = \dot{y}\tau + C_a$,

b. $\ddot{y}\tau^2 + 2\dot{y}\tau = 2y + C_b$,

c. $\ddot{y}\tau^2 + 2\dot{y}\tau = \dot{y}\tau + 2y + C_c$,

d. $\dot{y}\tau = 2y + C_d$,

e. $\dot{y}\tau = \dot{y}\tau + 2y + C_e$,

f. $2y = \dot{y}\tau + 2y + C_f$.

It is easy to show that every of condition a. – f. leads to trivial solution.

Let us consider Cauchy problem for the van der Pol equation solved as NLDE (the solution is (**)):
$$x_0 = x(t_0),$$
$$v_0 = \dot{x}(t_0),$$
determine $x(t)$ for all $t \geq t_0$. Recall that $\dot{x}(t_0) = \frac{x(t_0 + \tau) - x(t_0)}{\tau}$. From here and (**) one could find
$$v_0 = \frac{\frac{2+\lambda\tau}{\lambda\tau} - x_0^2}{x_0 \tau} = \frac{2 + \lambda\tau - x_0^2 \lambda\tau}{x_0 \lambda\tau^2}$$

and

$$\tau_\pm = \frac{\frac{1}{x_0} - x_0 \pm \sqrt{\left(x_0 - \frac{1}{x_0}\right)^2 + \frac{8v_0}{\lambda x_0}}}{2v_0},$$

that is $\tau$ is depend on the initial conditions. In linear cases (see examples above) it is not: the initial conditions determine integration constants values but $\tau$ is arbitrary there. The solution of Cauchy problem is one of the following:

$$x_{++} = \sqrt{\frac{2}{\lambda\tau_+} + 1} \left[\Omega_+ + \sqrt{\Omega_+^2 - 1}\right]^{\frac{(-1)^{t/\tau_+}}{2}},$$

$$x_{+-} = \sqrt{\frac{2}{\lambda\tau_+} + 1} \left[\Omega_+ - \sqrt{\Omega_+^2 - 1}\right]^{\frac{(-1)^{t/\tau_+}}{2}}, \quad (***)$$

$$x_{-+} = \sqrt{\frac{2}{\lambda\tau_-} + 1} \left[\Omega_- + \sqrt{\Omega_-^2 - 1}\right]^{\frac{(-1)^{t/\tau_-}}{2}},$$

$$x_{--} = \sqrt{\frac{2}{\lambda \tau_-} + 1} \left[ \Omega_- - \sqrt{\Omega_-^2 - 1} \right]^{\frac{(-1)^{t/\tau_-}}{2}}$$

where $\Omega_. = 1 + \frac{\omega^2 \tau_.^2}{2(2 + \lambda \tau_.)}$ and $\tau_\pm$ just defined. In any case we have:

$$x_\pm(t) = x_0^{(-1)^{\frac{(t-t_0)}{\tau_\pm}}},$$

where the explicit $\tau$ value is determined by initial conditions.

As a result, for any initial condition excluding $x_0 = 0$, $x_0 = 1$ and $v_0 = 0$, we have the only solution of the Cauchy problem which has one of the (***) form. This solution describes a chaotic dynamical system because of it obeys the following well-known criteria of the chaotic system: sensitivity to initial conditions, the presence of a topological mixing and everywhere dense periodic orbits. (The last two criteria may be inspected by construction of phase trajectories of the dynamical system.)

So, as we might expect, the discrete analysis of van der Pol oscillator leads to a consistent description and gives the exact solution of the Cauchy problem.

**Corollary**

The fundamental difference between the considered solutions of linear equations with one hand, and non-linear equation, on the other hand, lies in the fact that a specific value of non-limit differential step ($\tau$) does not have any influence on the form of solution in conventional (limited) form: for any step value in limitation $\tau \to 0$ will be the same expression. For non-linear equations (as van der Pol equation shows) it is not true in general. Namely:

*The nature of the equation, which describes the non-linear process, in general, depends on the length of time step in the process, that is, the value of $\tau$ (non-limit differential step), which is determined by the rest energy of the physical system together with the IRF of observer.*

Therefore, in general, any nonlinear physical process in a specific IRF has its own value of $\tau$, which is "encoded" by the equation of the process (equation of NLDE-type). Mathematically it means that if the independent variable is continuous, then the general solution sometimes can not be found in quadratures but partial solutions can be found. It is shown in case of van der Pol equation.

Physically this proposition means that laws are changeable with respect of different mass scales: micro, macro, mega.